\documentclass[10pt,a4paper]{article}
\usepackage{amsmath}\usepackage{theorem}\usepackage[all]{xy}
\usepackage{amssymb}\usepackage{latexsym}

\newtheorem{The}{Theorem}
\newtheorem{Pro}[The]{Proposition}

\newcommand{\bevis}{\textit{Proof: }}
\newcommand{\qed}{\hfill$\Box$}
\newcommand{\ner}{\vspace{\theorempostskipamount}}
\newcommand{\ra}{\rightarrow}
\newcommand{\la}{\leftarrow}

\newcommand{\Z}{\mathbf{Z}}
\newcommand{\N}{\mathbf{N}}
\newcommand{\pnt}{{\{\mathit{pt.}\}}}
\newcommand{\Top}{\mathbf{Top}}
\newcommand{\stj}{\,{\widetilde{*}}}

\begin{document}
\title{Bimonoidal operad-actions \\ and
the product in \\ negative Tate-cohomology  \ner\ner\ner
}

\author{\textsc{Pelle Salomonsson}}\date{\small{\texttt{pellegestalten@gmail.com}}}
\vspace{.1cm}
\maketitle
\vspace{.5cm}
\abstract{}
\!\!\!\!\!\!\!\noindent
\begin{quote} We study a certain construction designed to bring together the following two topics: $i$)~Dyer--Lashof-operations in negative Tate-cohomology, $ii$)~the description of negative Tate-cohomology in terms of joins. It has the merit of making (some) sense in a more general context: where the loop-space of the space under study no longer has to be of compact homotopy-type. The exposition given here is a streamlined version of a previous version of mine, available here at the Arxiv.
\end{quote}

\vspace{1cm}

\noindent The word ``bimonoidal'' in the title refers to a certain ``monoid/monoid-situation,'' not to a ``monoid/comonoid-situation.'' The latter sort of bimonoidalness is however (for seemingly unrelated reasons) present too (see section~\ref{section2}).

\section{The direct product and the join}
The monoid/monoid-situation occurs in the category of spaces, where we, apart from the disjoint union, have the two symmetric monoidal structures~$\times$ and~$*$, the direct product and the join. There is some distributivity between these (which we do not investigate in detail here, since we do not need it in its full generality).   More precisely, there are ``distributors''
\begin{gather}
(A*B)\times C\,\,\ra\,\,(A\times C)*(B\times C)  \,\,\,,  \label{distrib1}\\
(A\times B)* C\,\,\ra\,\,(A* C)\times(B*C)  \,\,\,,  \label{distrib2}
\end{gather}
both of them incidentally going in the same direction (the ``rewrite-direc\-tion''). Only the former one (the one that doesn't stabilise) can be of any use to us here, since it can be composed with the projection $A\times C\ra A$. We obtain a functorial ``priority-swapping-map''
\begin{gather}
(A*B)\times C\,\,\ra \,\,A*(B\times C) \,\,\,,
\end{gather}
whose existence, perhaps unlike~(\ref{distrib1}) and~(\ref{distrib2}), is quite clear to see. And it is precisely this map that is needed to give meaning the notion of a ``bimonoidal'' operad-action $\mathcal{C}(n)\times A^{*n}\ra A$. More precisely, it is needed to construct the shuffle-map in the pentagonal  diagram that appears among the axioms for an operad-action.

We fix a base-space~$X$, and consider the category of spaces mapping to it. There is a certain ``fibration-theoretic'' point of view, which often seems preferable. So let us restrict attention to \emph{fibrations} over~$X$, and let $\mathcal{C}$ be a $\Sigma$-free contractible operad in the category of spaces (as usual with respect to the direct product as monoidal structure). Then we consider the seemingly silly bimonoidal action $\mathcal{C}(n)\times\pnt^{*n} \ra \pnt$, where of course we might just as well  have written $\Delta^{n-1}$ for $\pnt^{*n}$\!. Given a space~$F$, we consider it as a ``homology-class'' in $\pnt$, and mimicking the construction of the Dyer--Lashof-action on infinite-loop-spaces, we obtain a mapping
\begin{equation}\label{fibrejoin}
\mathcal{C}(n)\times F^{*n} \,\,\ra\,\,\pnt\,\,\,.
\end{equation}
This may seem banal, but in our case this is a \emph{fibrewise} construction. Making everything vary in a family, the target-space is no longer trivial, and as the twisting is also non-trivial, we have an interesting functor that, upon passing to the quotient by a cyclic group $C_p\subset\Sigma_p$, to any $A\ra X$, a fibration say with fibre~$F$, associates a new space, denoted by $PA$ say, mapping to~$X$ with fibre $\mathcal{C}(p)\times_{C_p}F^{*p} $. This spatial construct may surely be of some independent interest, and at the level of homology we have:

\begin{Pro} \label{DL-proposition} Assume~$X$ to be the classifying-space of a finite group. If $A\ra X$ is the carrier of a distinguished mod~$p$ homo\-logy-class, then there is a naturally defined element in $H_{**}(PA)$, such that the association $[A]\mapsto[PA]$ defines an action of the mod~$p$ Dyer--Lashof-algebra on $H_*X$.
\end{Pro}
\bevis For the proof of this fact, one first of all needs the nonfibrewise product~$\stj$, to be defined immediately below. Then the details are rather clear, using ultimately the fact~\cite{tene} that Tate-cohomology can be defined in terms of the join. For example, the validity of the Adem-relations, in the mod~$2$ case say, follows ``as usual'' from the equivariant mapping
$\mathcal{C}(2)\times\mathcal{C}(2)\times\mathcal{C}(2)\times (F^{*2}*F^{*2})\ra \mathcal{C}(4)\times F^{*4}$.
\qed

\ner\noindent Apart from this fibrewise construction (of the join), there is another one: the union of paths going (as it were) from a point in~$A$ to one in~$B$ through~$X$, with the usual endpoint-identifications. More precisely, consider the space of points on paths $I\ra X$ together with the choice  of a point in~$A$ (and~$B$ respectively) mapping to each endpoint, and then perform the usual identifications (forgetting, at the~$A$-endpoint, about the rest of the path as well as the chosen point in~$B$,
and the other way around).

We may perhaps temporarily use the symbol~$\stj$ to denote this other construction. But most of the time, we would surely prefer to work up to homotopy, and leave it unspecified which construction is intended. Both are clearly homotopy-invariant, taken separately. Moreover, we have the following comparison.
\begin{Pro} For  two fibrations $A\ra X\la B$, the fibrewise join  $A*_XB$ is homotopy-equivalent to the nonfibrewise one $A\stj_XB$.
\end{Pro}
\bevis To see that the inclusion
$A*_XB\ra A\stj_XB$  is a homotopy-equival\-ence, we must exhibit a scheme to systematically  move general paths to constant ones. To achieve that, let us (by homotopy-invariance) assume both spaces to be given to us as pathspace-fibrations. That is, we assume that $A=PA'$ and $B=PB'$. Then the required moving-strategy is readily found. Notice: the usage of the symbol ``$P$'' here is of course unrelated to the one in proposition~\ref{DL-proposition}. \qed

\ner\noindent The importance of the nonfibrewise construction lies in the existence of yet another product, which cannot be defined in the fibrewise manner (should we try to do that, we would just, up to homotopy, get  the fibrewise direct product). It is obtained from~$\stj$ by forgetting slightly less at the endpoints. Namely, at the $A$-endpoint we remember the chosen point in~$B$, and conversely, although we still forget about the path. We may perhaps use the symbol~$\bullet$ to denote this product, which is quite ``bad'' in the sense that it is nonassociative and impossible to define in a fibrewise way. It is however of great importance for the definition of the product in negative (field-coefficient) Tate-cohomology. Namely, we have the diagram
\begin{equation}\label{birational}
A\stj_X B\,\,\la \,\,A\bullet_X B\,\,\ra \,\,A\times B\,\,\,,
\end{equation}
where the unindexed symbol~$\times$ refers to the absolute product (the product in the category of spaces), not the relative version. The importance of that diagram lies in certain good features that the mappings in it have,  making them suitable for homology-considerations. Indeed, if~$A$ and~$B$ are finite cell-complexes, and if we are given a finitedimensional model of $\Omega X$, then the left-hand map restricts to a homeomorphism over a dense open subset. And as for the right-hand map, it is a fibration with fibre $\Sigma\Omega X$.

\section{The product and the coproduct in negative Tate-cohomology}\label{section2}
We may perhaps proceed to say a few words about how the cup-product, and hence ``Tate-duality,'' may be fitted into this description. We end up with the diagram
\begin{equation}\xymatrix@R+.3pc@C+1.3pc{
X\times X& \txt{$X\times X$\\$\amalg$\\$X\times X$}
\ar[r]_-{\mathit{diag.}~\times~\mathit{id.}}
^-{\mathit{id.}~\times~\mathit{diag.}}\ar[l]_{\mathit{collapse}}&
\txt{$X\times X\times X$\\$\amalg$\\$X\times X\times X$}\\
X\bullet_XX\ar@{->}[u]\ar[d]&  &
\txt{$(X\bullet_XX)\times X$\\$\amalg$\\$X\times (X\bullet_XX)$}\ar@{->}[u]\ar[d]\\
X\ar[rr]& &X\times X
}  \label{frobeniusdiagram}
\end{equation}
which admittedly may be accused of a certain awkwardness. Some of the maps go in the ``wrong'' direction, and the $\bullet$-operation is nonassociative. And not only are there \emph{two} different binary operations present (apart from the disjoint union), but they seem much too unrelated to each other, being defined in two different categories ($\Top$ and $\Top_X$ respectively). Some of these deficiencies may certainly be remedied by recasting the diagram into another form, but in the final end their presence probably indicates that there cannot be any genuine ``Frobenius-objects'' in the category of spaces (in contrast thus to the situation concerning ``Hopf-objects''). We have however the following fact, which we state somewhat below its natural level of generality (where the finite group should be replaced by a compact Lie-group).
\begin{Pro} Assume~$X$ to be the classifying space of a finite group. Then diagram~(\ref{frobeniusdiagram}) commutes upon passing to homology (with coefficients in a field). That is to say: it commutes when one begins in the upper left-hand corner, and ends in the lower right-hand one (and upon using wrong-way-maps when necessary).
\end{Pro}
\bevis We use the known fact that Tate-cohomology is a Frobenius-algebra. Indeed, this is one way to phrase the existence of ``Tate-duality.'' It follows from the formalism of Frobenius-algebras that we get a diagram of the stated form, except that one would not expect the redoublings (the disjoint unions) in the exhibited diagram to be there. They are however caused by the passing from the $\Z$-graded structure $\widehat{H}^*X$ to the $\N$-graded one $H_*X$. More precisely, beginning with $a\otimes b$ in the homology of the upper left-hand corner, and using lower indices to write graded components of the coproduct (assuming them to be tensor-monomials), we obtain along the upper path
$\sum ab_{n'}\otimes b_{n''}+\sum a_{m'}\otimes a_{m''}b$, where $m:=m'+m''$ and $n:=n'+n''$ are the degrees of~$a$ and~$b$, and where juxtaposition denotes multiplication in negative Tate-cohomology (which here increases degree by one). The two summands are in fact equal, and they are moreover, each one of them, equal to the result along the lower path. Hence there seems to be a discrepancy by a factor of~$2$. However, for each bidegree only one (exactly one) of the summands survives into the $\N$-graded structure. Namely, equality of bidegrees entails that $n'+m''=-1$. Hence exactly one of them is $\geq0$. \qed

\ner\ner\noindent\textbf{Remark. Overall motivation.} As I said in the abstract, these developments are designed to bring together the following two topics: $i$)~Dyer--Lashof-operations in negative Tate-cohomology, $ii$)~the description of negative Tate-cohomology in terms of joins. As for ($i$), see Langer's paper~\cite{langer}. As for ($ii$), see Haggai Tene's  paper~\cite{tene}. The recourse in~\cite{salomonsson} to a particular operad-model was the result of confusion.

\vspace{1cm}


\end{document}